# Integral Expressions for $\zeta_2(s_1,s_2)$, $\zeta_3(s_1,s_2,s_3)$ and Approximations

VIVEK V. RANE

A-3/203, ANAND NAGAR,

DAHISAR, MUMBAI-400 068

INDIA.

v_v_rane@yahoo.co.in

**Abstract :** For the Euler-Zagier's multiple zeta function $\zeta_2(s_1,s_2)$ of two variables, we obtain its integral representation (involving product of Hurwitz's zeta functions) over the interval $[1,\infty)$, with respect to the second variable of the Hurwitz zeta function and also obtain a good approximation to it as a function of $s_1$ and $s_2$, for $s_1 \geq 1$ and $s_2 > 1$. In particular, this approximate is explicitly computable, when $s_1$, $s_2$ differ by an even integer and is good, when $s_2$ is large. We treat $\zeta_3(s_1,s_2,s_3)$ likewise.

**Keywords :** Hurwitz's /Euler-Zagier's/Tornheim /Riemann zeta functions.

# Integral Expressions for $\zeta_2(s_1,s_2)$, $\zeta_3(s_1,s_2,s_3)$ and Approximations

VIVEK V. RANE,

A-3/203, ANAND NAGAR,

DAHISAR, MUMBAI-400 068

INDIA.

v_v_rane@yahoo.co.in

For the complex variable $s$ and a fixed complex number $\alpha \neq 0,-1,-2,\ldots$ , define the Hurwitz's zeta function $\zeta(s,\alpha)$ by $\zeta(s,\alpha) = \sum_{n \geq 0}(n+\alpha)^{-s}$ for Re s>1; and its analytic continuation as a function of $s$ and treat $\alpha$ as a complex variable thereafter. We shall write $\zeta(s,1) = \zeta(s)$, the Riemann zeta function. For the complex variable $s$ with Re s>1 and for the real variable $\alpha$, we define $\zeta^{(r)}(s,\infty) = \lim_{\alpha \to \infty}\zeta^{(r)}(s,\alpha)$, where $\zeta^{(r)}(s,\alpha) = \dfrac{\partial^r}{\partial s^r}\zeta(s,\alpha)$ for any integer $r \geq 0$. Note that for Re s>1, $\zeta^{(r)}(s,\infty) = 0$. In particular, $\zeta(s,\infty) = 0$ for Re s>1. Also note that for s>1, $\zeta(s,u)$ is a monotonically decreasing function of $u$.

For the complex variables $s_1, s_2, \ldots, s_r$, where $r \geq 1$ is an integer, define the multiple zeta function $\zeta_r(s_1,s_2\ldots s_r) = \sum_{n_1 \geq 1}n_1^{-s_1}\sum_{n_2>n_1}n_2^{-s_2}\ldots\sum_{n_r>n_{r-1}}n_r^{-s_r}$

for Re $s_i \geq 1$ for $i = 1,2,\ldots,r-1$ and Re $s_r > 1$; and its analytic continuation . In

particular $\zeta_2(s_1, s_2) = \sum_{n_1 \geq 1} n_1^{-s_1} \sum_{n_2 > n_1} n_2^{-s_2}$ for Re $s_1 \geq 1$ and Re $s_2 > 1$; and its analytic

continuation .

For the complex variables $s_1, s_2, s_3$, Tornheim's double zeta function

T($s_1, s_2, s_3$) is defined by $T(s_1, s_2, s_3) = \sum_{n_1 \geq 1} \sum_{n_2 \geq 1} n_1^{-s_1} n_2^{-s_2} (n_1 + n_2)^{-s_3}$ for Re($s_1 + s_3$) > 1,

Re($s_2 + s_3$) > 1 and Re($s_1 + s_2 + s_3$) > 2 ; and its analytic continuation . Then T($s_1, 0, s_3$)

= $\zeta_2$ ($s_1, s_3$) . Using the theory of Tornheim's double zeta function , one can evaluate

T($n_1, 0, n_2$) for the integers $n_1, n_2 \geq 2$ with their difference an odd integer .

However, for the integers $n_1, n_2 \geq 2$ with an even difference , the evaluation is

complicated . See for example author [2] . In this backdrop , we shall give an

integral representation for $\zeta_2(s_1, s_2)$ for $s_1, s_2 > 1$ and give an approximation to it in

terms of integrals (involving products of Hurwitz zeta functions) over the interval

$[1, \infty)$ with respect to the second variable .Interestingly , this approximate is

explicitly computable, when $s_1 \geq 1, s_2 > 1$ differ by an even integer and is very good for large $s_2$. Note that for any complex $s$ with Re $s>1$, the integral

$$\int_1^\infty \zeta(s,u)\zeta(s+1,u)du = -\frac{1}{2s}\int_1^\infty \frac{d}{du}\zeta^2(s,u)du = -\frac{1}{2s}\left(\zeta^2(s,\infty) - \zeta^2(s,1)\right) = \frac{\zeta^2(s)}{2s}.$$

In view of this, on integration by parts repeatedly, the integral $\int_1^\infty \zeta(s,u)\zeta(s+n,u)du$ can be explicitly evaluated in terms of products of values of $\zeta(s)$ for Re $s>1$, when $n$ is an odd positive integer. Similarly, in the context of $\zeta_3(s_1, s_2, s_3)$, it is to be noted that for $s>1$,

$$\int_1^\infty \zeta^2(s,u)\zeta(s+1,u)du = -\frac{1}{3s}\int_1^\infty \frac{d}{du}\zeta^3(s,u)du = -\frac{1}{3s}\left(\zeta^3(s,\infty) - \zeta^3(s,1)\right) = \frac{\zeta^3(s)}{3s}.$$

See for example author [1], where we have developed the integral calculus of $\frac{\partial^r}{\partial s^r}\zeta(s,\alpha)$ with respect to $\alpha$. Against this overall backdrop, we shall state our Theorem 1 and Theorem 2 below.

In what follows, [u] shall denote the integral part of the real variable u and $\zeta(s)$ shall stand for Riemann zeta function.

**Theorem 1 :** For the real variables $s_1 \geq 1$, $s_2 > 1$, we have

$$\zeta_2(s_1, s_2) = s_1 \int_1^\infty \zeta(s_1+1, u)\zeta(s_2, [u]+1)du = s_1\int_1^\infty \zeta(s_1+1, u)\zeta(s_2, u)du - A(s_1, s_2),$$

where $0 \leq A(s_1, s_2) \leq s_1 \int_1^\infty \frac{\zeta(s_1+1, u) du}{u^{s_2}} \leq \frac{s_1}{s_2 - 1} \zeta(s_1 + 1)$

**Corollary :** We have $\zeta_2(1, s_2) = \int_1^\infty \zeta(s_2, u) \zeta(2, [u]+1) du = \int_1^\infty \zeta(s_2, u) \zeta(2, u) du - A(1, s_2)$ ,

where $0 \leq A(1, s_2) \leq \frac{\zeta(2)}{s_2 - 1}$ .

**Theorem 2 :** Let $s_1, s_2, s_3 > 1$ . Then we have

$$\zeta_3(s_1, s_2, s_3) = \zeta(s_1)\zeta(s_2)\zeta(s_3) - \zeta_2(s_1, s_2 + s_3)$$

$$- \zeta_2(s_1 + s_2, s_3) - \zeta(s_1)\zeta_2(s_3, s_2) - \zeta(s_1 + s_2 + s_3) - s_2 \int_1^\infty \zeta(s_1, [u]+1)\zeta(s_2+1, u)\zeta(s_3, u) du$$

$$- s_3 \int_1^\infty \zeta(s_1, [u]+1)\zeta(s_2, u)\zeta(s_3+1, u) du + s_3 \int_1^\infty \zeta(s_1, [u]+1)\zeta(s_2, [u]+1)\zeta(s_3+1, u) du .$$

**Proof of Theorem 1 :** We have for $s_1, s_2 > 1$ ,

$$\zeta_2(s_1, s_2) = \sum_{n_1 \geq 1} n_1^{-s} \cdot \sum_{n_2 > n_1} n_2^{-s_2} = \sum_{n_1 \geq 1} n_1^{-s_1} \zeta(s_2, n_1 + 1)$$

$$= \sum_{n_1 \geq 1} n_1^{-s_1} \left( \zeta(s_2, n_1) - n_1^{-s_2} \right) = \sum_{n \geq 1} n^{-s_1} \zeta(s_2, n) - \sum_{n \geq 1} n^{-(s_1+s_2)} = \sum_{n \geq 1} n^{-s} \zeta(s_2, n) - \zeta(s_1 + s_2).$$

Next , we shall evaluate $\sum_{n \geq 1} n^{-s_1} \zeta(s_2, n) = \zeta(s_2) + \sum_{n \geq 2} n^{-s_1} \zeta(s_2, n)$ .

Next $\sum_{n\geq 2} n^{-s_1}\zeta(s_2,n) = -\sum_{n\geq 2}\zeta(s_2,n)(\zeta(s_1,n+1)-\zeta(s_1,n))$ .

Note that $\zeta(s_1,n+1)-\zeta(s_1,n)$ is the jump in value of the step function $\zeta(s_1,[u]+1)$

at the point u=n for any integer $n\geq 2$ . In view of the definition of

Riemann-Stieltjes integral of a continuous function with respect to a step function as integrator , we have $\sum_{n\geq 2}\zeta(s_2,n)(\zeta(s_1,n+1)-\zeta(s_1,n)) = \int_1^\infty \zeta(s_2,u)d\zeta(s_1,[u]+1)$ ,

where we have defined $\int_1^\infty \zeta(s_2,u)d\zeta(s_1,[u]+1) = \lim_{N\to\infty}\int_1^N \zeta(s_2,u)d\zeta(s_1,[u]+1)$ ,

after noting that $\lim_{u\to\infty}\zeta(s,u)=0$ for s>1 .

Thus $\sum_{n\geq 2} n^{-s_1}\zeta(s_2,n)$

$$= -\int_1^\infty \zeta(s_2,u)d\zeta(s_1,[u]+1) = -\left\{\zeta(s_2,u)\zeta(s_1,[u]+1)\Big|_{u=1}^\infty - \int_1^\infty \zeta(s_1,[u]+1)d\zeta(s_2,u)\right\}$$

$$= -\left\{\zeta(s_2,\infty)\zeta(s_1,\infty) - \zeta(s_2)\zeta(s_1,2) + s_2\int_1^\infty \zeta(s_1,[u]+1)\zeta(s_2+1,u)\right\}$$

$$= \zeta(s_2)(\zeta(s_1)-1) - s_2\int_1^\infty \zeta(s_1,[u]+1)\zeta(s_2+1,u)du$$

$$= \zeta(s_1)\zeta(s_2) - \zeta(s_2) - s_2 \int_1^\infty \zeta(s_1,[u]+1)\zeta(s_2+1,u)du .$$

Thus $\sum_{n\geq 1} n^{-s_1}\zeta(s_2,n) = \zeta(s_1)\zeta(s_2) - s_2 \int_1^\infty \zeta(s_1,[u]+1)\zeta(s_2+1,u)du .$

Thus $\zeta_2(s_1,s_2) = \zeta(s_1)\zeta(s_2) - \zeta(s_1+s_2) - s_2 \int_1^\infty \zeta(s_1,[u]+1)\zeta(s_2+1,u)du .$

This gives on interchanging $s_1$, $s_2$,

$$\zeta_2(s_2,s_1) = \zeta(s_1)\zeta(s_2) - \zeta(s_1+s_2) - s_1 \int_1^\infty \zeta(s_2,[u]+1)\zeta(s_1+1,u)du .$$

Next, we know that $\zeta_2(s_1,s_2) + \zeta_2(s_2,s_1) = \zeta(s_1)\zeta(s_2) - \zeta(s_1+s_2)$ .

This gives $\zeta_2(s_1,s_2) = s_1 \int_1^\infty \zeta(s_1+1,u)\zeta(s_2,[u]+1)du$ for $s_1,s_2>1$ .

Next note that $u-1 \leq [u] \leq u$ for any real $u$ so that $u \leq [u]+1 \leq u+1$ .

This gives for $s_2>1$ , $\zeta(s_2,u) \geq \zeta(s_2,[u]+1) \geq \zeta(s_2,u+1)$ .

This gives $s_1 \int_1^\infty \zeta(s_2,u)\zeta(s_1+1,u)du$

$$\geq s_1 \int_1^\infty \zeta(s_2,[u]+1)\zeta(s_1+1,u)du \geq s_1 \int_1^\infty \zeta(s_2,u+1)\zeta(s_1+1,u)du .$$

Note $\zeta(s_2, u+1) = \zeta(s_2, u) - u^{-s_2}$ .

This gives $s_1 \int_1^\infty \zeta(s_2, u)\zeta(s_1+1, u)du - s_1 \int_1^\infty \zeta(s_1, +1, u)u^{-s_2} du$

$\leq \zeta_2(s_1, s_2) \leq s_1 \int_1^\infty \zeta(s_2, u)\zeta(s_1+1, u)du$ .

This gives $\zeta_2(s_1, s_2) = s_1 \int_1^\infty \zeta(s_2, u)\zeta(s_1+1, u)du - A(s_1, s_2)$ ,

where $0 \leq A(s_1, s_2) \leq s_1 \int_1^\infty \zeta(s_1+1, u)u^{-s_2} du \leq s_1\zeta(s_1+1)\int_1^\infty \frac{du}{u^{s_2}} \leq \frac{s_1}{s_2-1}\zeta(s_1+1)$ .

In particular , letting $s_1 \to 1$ from the right , we get

$\zeta_2(1, s_2) = \int_1^\infty \zeta(s_2, u)\zeta(2, u)du - A(1, s_2)$, where $A(1, s_2) \leq \frac{1}{s_2-1}\zeta(2)$ .

**Proof of Theorem 2 :** We have $\zeta_3(s_1, s_2, s_3) = \sum_{n_1 \geq 1} n_1^{-s_1} \cdot \sum_{n_2 > n_1} n_2^{-s_2} \cdot \sum_{n_3 > n_2} n_3^{-s_3}$

Consider $\sum_{n_2 > n_1} n_2^{-s_2} \sum_{n_3 > n_2} n_3^{-s_3} = \sum_{n_2 > n_1} n_2^{-s_2} \zeta(s_3, n_2+1) = \sum_{n_2 > n_1} n_2^{-s_2} \left(\zeta(s_3, n_2) - n_2^{-s_3}\right)$

$= \sum_{n_2 > n_1} n_2^{-s_2} \zeta(s_3, n_2) - \sum_{n_2 > n_1} n_2^{-(s_2+s_3)}$ .

Thus $\zeta_3(s_1, s_2, s_3) = \sum_{n_1 \geq 1} n_1^{-s_1} \sum_{n_2 > n_1} n_2^{-s_2} \sum_{n_3 > n_2} n_3^{-s_3} = \sum_{n \geq 1} n_1^{-s_1} \sum_{n_2 > n_1} n_2^{-s_2} \cdot \zeta(s_3, n_2) - \sum_{n_1 \geq 1} n_1^{-s_1} \sum_{n_2 > n_1} n_2^{-(s_2+s_3)}$

$$= \sum_{n_1 \geq 1} n_1^{-s_1} \cdot \sum_{n_2 > n_1} n_2^{-s_2} \cdot \zeta(s_3, n_2) - \zeta_2(s_1, s_2 + s_3).$$

Next consider

$$\sum_{n_2 > n_1} n_2^{-s_2} \zeta(s_3, n_2) = -\sum_{n_2 > n_1} \zeta(s_3, n_2)\big(\zeta(s_2, n_2+1) - \zeta(s_2, n_2)\big) = -\int_{n_1}^{\infty} \zeta(s_3, u) d\zeta(s_2, [u]+1)$$

$$= -\left\{ \zeta(s_3, u)\zeta(s_2, [u]+1) \big|_{u=n_1}^{\infty} + s_3 \int_{n_1}^{\infty} \zeta(s_2, [u]+1)\zeta(s_3+1, u) du \right\}$$

$$= -\left\{ \zeta(s_3, \infty)\zeta(s_2, \infty) - \zeta(s_3, n_1)\zeta(s_2, n_1+1) + s_3 \int_{n_1}^{\infty} \zeta(s_2, [u]+1)\zeta(s_3+1, u) du \right\}$$

$$= \zeta(s_3, n_1)\zeta(s_2, n_1+1) - s_3 \int_{n_1}^{\infty} \zeta(s_2, [u]+1)\zeta(s_3+1, u) du.$$

Thus $\sum_{n_1 \geq 1} n_1^{-s_1} \sum_{n_2 > n_1} n_2^{-s_2} \cdot \zeta(s_3, n_2) = \sum_{n_1 \geq 1} n_1^{-s_1} \cdot \zeta(s_2, n_1)\zeta(s_3, n_1) - \sum_{n_1 \geq 1} n_1^{-(s_1+s_2)} \cdot \zeta(s_3, n_1)$

$$- s_3 \sum_{n_1 \geq 1} n_1^{-s_1} \cdot \int_{n_1}^{\infty} \zeta(s_2, [u]+1)\zeta(s_3+1, u) du$$

Next $\sum_{n_1 \geq 1} n_1^{-s_1} \cdot \zeta(s_2, n_1)\zeta(s_3, n_1) = \zeta(s_2)\zeta(s_3) + \sum_{n_1 \geq 2} \zeta(s_2, n_1)\zeta(s_3, n_1) n_1^{-s_1}$

$$= \zeta(s_2)\zeta(s_3) - \left\{ \sum_{n_1 \geq 2} \zeta(s_2, n_1)\zeta(s_3, n_1) \cdot \left(\zeta(s_1, n_1 + 1) - \zeta(s_1, n_1)\right) \right\}$$

$$= \zeta(s_2)\zeta(s_3) - \int_1^\infty \zeta(s_2, u)\zeta(s_3, u) d\zeta(s_1, [u]+1)$$

$$= \zeta(s_2)\zeta(s_3) - \left\{ \zeta(s_2, u)\zeta(s_3, u)\zeta(s_1, [u]+1) \Big|_{u=1}^\infty - \int_1^\infty \zeta(s_1, [u]+1) d(\zeta(s_2, u)\zeta(s_3, u)) \right\}$$

$$= \zeta(s_2)\zeta(s_3) - \left\{ -\zeta(s_2)\zeta(s_3)\zeta(s_1, 2) - \int_1^\infty \zeta(s_1, [u]+1) d(\zeta(s_2, u)\zeta(s_3, u)) \right\}$$

$$= \zeta(s_1)\zeta(s_2)\zeta(s_3) + \int_1^\infty \zeta(s_1, [u]+1) d(\zeta(s_2, u)\zeta(s_3, u))$$

Next $\sum_{n_1 \geq 1} n_1^{-(s_1+s_2)} \cdot \zeta(s_3, n_1) = \zeta_2(s_1+s_2, s_3) + \zeta(s_1+s_2+s_3)$.

Next

$$\sum_{n_1 \geq 1} n_1^{-s_1} \int_{n_1}^\infty \zeta(s_2, [u]+1)\zeta(s_3+1, u) du = \sum_{n_1 \geq 1} n_1^{-s_1} \left( \int_1^\infty \zeta(s_2, [u]+1)\zeta(s_3+1, u) du - \int_1^{n_1} \zeta(s_2, [u]+1)\zeta(s_3+1, u) du \right)$$

$$= \sum_{n_1 \geq 1} n_1^{-s_1} \cdot \int_1^\infty \zeta(s_2, [u]+1)\zeta(s_3+1, u) du - \sum_{n_1 \geq 1} n_1^{-s_1} \int_1^{n_1} \zeta(s_2, [u]+1)\zeta(s_3+1, u) du$$

$$= \zeta(s_1) \int_1^\infty \zeta(s_2, [u]+1)\zeta(s_3+1, u) du - \sum_{n_1 \geq 1} n_1^{-s_1} \int_1^{n_1} \zeta(s_2, [u]+1)\zeta(s_3+1, u) du$$

$$= \frac{\zeta(s_1)}{s_3} \cdot \zeta_2(s_3, s_2) - \sum_{n_1 \geq 1} n_1^{-s_1} \cdot \int_1^{n_1} \zeta(s_2, [u]+1)\zeta(s_3+1, u)du, \text{on using Theorem 1}.$$

Thus $-s_3 \sum_{n_1 \geq 1} n_1^{-s_1} \int_{n_1}^{\infty} \zeta(s_2, [u]+1)\zeta(s_3+1, u)du$

$$= s_3 \sum_{n_1 \geq 1} n_1^{-s_1} \cdot \int_1^{n_1} \zeta(s_2, [u]+1)\zeta(s_3+1, u)du - \zeta(s_1)\zeta_2(s_3, s_2)$$

This gives $\sum_{n_1 \geq 1} n_1^{-s_1} \sum_{n_2 > n_1} n_2^{-s_2} \zeta(s_3, n_2) = \zeta(s_1)\zeta(s_2)\zeta(s_3)$

$$+ \int_1^{\infty} \zeta(s_1, [u]+1) d(\zeta(s_2, u)\zeta(s_3, u)) - \zeta_2(s_1+s_2, s_3) - \zeta(s_1+s_2+s_3)$$

$$+ s_3 \sum_{n_1 \geq 1} n_1^{-s_1} \cdot \int_1^{n_1} \zeta(s_2, [u]+1) \cdot \zeta(s_3+1, u)du - \zeta(s_1)\zeta_2(s_3, s_2)$$

so that $\zeta_3(s_1, s_2, s_3) = \zeta(s_1)\zeta(s_2)\zeta(s_3) - \zeta_2(s_1, s_2+s_3)$

$$- \zeta_2(s_1+s_2, s_3) - \zeta(s_1)\zeta_2(s_3, s_2) - \zeta(s_1+s_2+s_3)$$

$$+ \int_1^{\infty} \zeta(s_1, [u]+1) d(\zeta(s_2, u)\zeta(s_3, u)) + s_3 \sum_{n_1 \geq 1} n_1^{-s_1} \int_1^{n_1} \zeta(s_2, [u]+1)\zeta(s_3+1, u)du$$

Next $\sum_{n_1 \geq 1} n_1^{-s_1} \cdot \int_1^{n_1} \zeta(s_2, [u]+1)\zeta(s_3+1, u)du = -\sum_{n_1 \geq 2} (\zeta(s_1, n_1+1) - \zeta(s_1, n_1)) \cdot \int_1^{n_1} \zeta(s_2, [v]+1)\zeta(s_3+1, v)dv$

$$= -\int_1^\infty (\int_1^u \zeta(s_2,[v]+1)\zeta(s_3+1,v)dv)d\zeta(s_1,[u]+1)$$

$$= -\left\{ \zeta(s_1,[u]+1)\cdot \int_1^u \zeta(s_2,[v]+1)\zeta(s_3+1,v)dv \Big|_{u=1}^\infty - \int_1^\infty \zeta(s_1,[u]+1)\zeta(s_2,[u]+1)\zeta(s_3+1,u)du \right\}$$

$$= -\left\{ \zeta(s_1,\infty)\int_1^\infty \zeta(s_2,[v]+1)\zeta(s_3+1,v)dv - \zeta(s_1,2)\cdot \int_1^1 \zeta(s_2,[v]+1)\zeta(s_3+1,v)dv \right.$$

$$\left. - \int_1^\infty \zeta(s_1,[u]+1)\zeta(s_2,[u]+1)\zeta(s_3+1,u)du \right\}$$

$$= \int_1^\infty \zeta(s_1,[u]+1)\zeta(s_2,[u]+1)\zeta(s_3+1,u)du$$

Thus $\zeta_3(s_1,s_2,s_3) = \zeta(s_1)\zeta(s_2)\zeta(s_3) - \zeta_2(s_1,s_2+s_3)$

$$-\zeta_2(s_1+s_2,s_3) - \zeta(s_1)\zeta_2(s_3,s_2) - \zeta(s_1+s_2+s_3)$$

$$- s_2\int_1^\infty \zeta(s_1,[u]+1)\zeta(s_2+1,u)\zeta(s_3,u)du - s_3\int_1^\infty \zeta(s_1,[u]+1)\zeta(s_2,u)\zeta(s_3+1,u)du$$

$$+ s_3\int_1^\infty \zeta(s_1,[u]+1)\zeta(s_2,[u]+1)\zeta(s_3+1,u)du$$

This completes the proof of Theorem 2.